\documentclass[11pt]{article}
\usepackage{amsmath, amsthm}

\def\tr{\mathop{\rm tr }\nolimits}

\def\qed{{\ }\hfill$\rlap{$\sqcup$}\sqcap$\smallskip}

\def\1{{\bf 1}}
\def\0{{\bf 0}}
\def\2{{\textstyle \frac{1}{2}}}

\def\R{{\rm I\kern-.21em R}}
\newtheorem*{thm}{Theorem}
\newtheorem*{lem}{Lemma}
\newtheorem*{prp}{Proposition}
\newtheorem*{cor}{Corollary}
\newtheorem*{sethm}{Spectral Excess Theorem}

\begin{document}
\author{Edwin R. van Dam
and Willem H. Haemers\\
{\footnotesize Tilburg University, Dept. Econometrics and O.R.,}\\[-3pt]
{\footnotesize P.O. Box 90153, 5000~LE~~Tilburg, The Netherlands,}\\[-3pt]
{\footnotesize e-mail: Edwin.vanDam@uvt.nl, Haemers@uvt.nl}
}
\title{An odd characterization of the generalized odd graphs\footnote{This version is published in Journal of Combinatorial Theory, Series B 101 (2011), 486-489.}}
\date{}
\maketitle \noindent {\footnotesize {\bf Abstract} We show that any
connected regular graph with $d+1$ distinct eigenvalues and odd-girth $2d+1$
is distance-regular, and in particular that it is a generalized odd graph.\\


}

\section{Introduction}\label{i}

The odd-girth of a graph is the length of the shortest odd cycle. A generalized
odd graph is a distance-regular graph of diameter $D$ and odd girth $2D+1$. It
is also called an almost-bipartite distance-regular graph, or a regular thin
near $(2D+1)$-gon. Well-known examples of such graphs are the Odd graphs (also
known as the Kneser graphs $K(2D+1,D)$), and the folded $(2D+1)$-cubes.

In this note, we shall characterize these graphs, by showing that any connected
regular graph with $d+1$ distinct eigenvalues and odd-girth (at least) $2d+1$
is a distance-regular generalized odd graph. We remark that $D=d$ for
distance-regular graphs, but for arbitrary connected graphs we only have the
inequality $D \le d$. In general it is not true that any connected regular
graph with diameter $D$ and odd-girth $2D+1$ is a generalized odd graph.
Counterexamples can easily be found for the case $D=2$ (among the triangle-free
regular graphs with diameter two there are many graphs that are not strongly
regular).

Huang and Liu \cite{HL99} proved that any graph with the same
spectrum as a generalized odd graph is such a graph. Because
the odd-girth of a graph follows from the spectrum, our
characterization is a generalization of this result.

For background on distance-regular graphs we refer the reader
to \cite{bcn}, for eigenvalues of graphs to \cite{CDS}, for
spectral characterizations of graphs to \cite{DH, DH09}, and
for spectral and other algebraic characterizations of
distance-regular graphs to \cite{DHKS} and \cite{Fi02},
respectively. To show the claimed characterization, we shall
use the so-called spectral excess theorem due to Fiol and
Garriga \cite{FG97}. Let $\Gamma$ be a connected $k$-regular
graph with $d+1$ distinct eigenvalues. The {\em excess} of a
vertex $u$ of $\Gamma$ is the number of vertices at distance
$d$ from $u$. We also need the so-called predistance polynomial
$p_d$ of $\Gamma$, which will be explained in some detail in
Section~\ref{odd-girth}. The important property of $p_d$ is
that the value of $p_d(k)$ --- the so-called {\em spectral
excess} --- only depends on the spectrum of $\Gamma$ (in fact,
all predistance polynomials depend only on the spectrum).
\begin{sethm} 
Let $\Gamma$ be a connected regular graph with $d+1$ distinct
eigenvalues. Then $\Gamma$ is distance-regular if and only if the
average excess equals the spectral excess.
\end{sethm}
\noindent For short proofs of this theorem we refer the reader
to \cite{Damexcess, Fiolexcess}. Note that one can even show
that the average excess is at most the spectral excess, and
that in \cite{Damexcess}, a bit stronger result is obtained by
using the harmonic mean of the number of vertices minus the
excess, instead of the arithmetic mean.

\section{The spectral characterization}\label{HL}

Let $\Gamma$ be a connected $k$-regular graph with adjacency
matrix $A$ having $d+1$ distinct eigenvalues
$k=\lambda_0>\lambda_1> \cdots > \lambda_d$ and finite
odd-girth at least $2d+1$. It follows that every vertex $u$ has
vertices at distance $d$, because otherwise the vertices at odd
distance from $u$ on one hand and the vertices at even distance
from $u$ on the other hand, would give a bipartition of the
graph, contradicting that the odd-girth is finite. Because
$\Gamma$ has diameter $D$ at most $d$, it follows that $D=d$,
and that the odd-girth equals $2d+1$.

Because $(A^i)_{uv}$ counts the number of walks of length $i$
in $\Gamma$ from $u$ to $v$, it follows that $p(A)$ has zero
diagonal for any odd polynomial $p$ of degree at most $2d-1$.
Therefore also $\tr p(A) = 0$. Because the trace of $p(A)$ can
also be expressed in terms of the spectrum of $\Gamma$, this
also shows that the odd-girth condition on $\Gamma$ is a
condition on the spectrum of $\Gamma$. In the following, we
make frequent use of polynomials. One of these is the Hoffman
polynomial $H$ defined by $H(x)=\frac n{\pi_0}\prod_{i=1}^d
(x-\lambda_i)$, where $n$ is the number of vertices and
$\pi_0=\prod_{i=1}^d (k-\lambda_i)$. This polynomial satisfies
$H(A)=J$, the all-ones matrix.

Let us now consider two arbitrary vertices $u,v$ at distance
$d$. By considering the Hoffman polynomial, it follows that
$(A^{d})_{uv}=\frac{\pi_0}n$.
By considering the minimal polynomial (or $(x-k)H$), it follows
that $(A^{d+1})_{uv}-\tilde{a}_d(A^d)_{uv}=0$, where
$\tilde{a}_d=\sum_{i=0}^d \lambda_i$ is the coefficient of
$x^d$ in the minimal polynomial.
Hence $(A^{d+1})_{uv}=\tilde{a}_d\frac{\pi_0}n$.
\begin{lem}
The average excess $\overline{k_d}$ of $\Gamma$ equals
$\frac{n}{\tilde{a}_d \pi_0^2}\tr A^{2d+1}$.
\end{lem}
\noindent {\bf Proof.} For a vertex $u$, let $\Gamma_d(u)$ be the
set of vertices at distance $d$ from $u$. Then
$$(A^{2d+1})_{uu}=\sum_{v \in \Gamma_d(u)}(A^d)_{uv}(A^{d+1})_{vu}
= k_d(u)\tilde{a}_d \pi_0^2/{n^2},$$ where $k_d(u) = |\Gamma_d(u)|$
is the excess of $u$. Therefore $\overline{k_d}\tilde{a}_d\pi_0^2/n
= \tr A^{2d+1}$ and $\tilde{a}_d \neq 0$.~\qed
\\[5pt]
In order to apply the spectral excess theorem, we have to
ensure that $\overline{k_d} = p_d(k)$. However, $p_d(k)$ and
$\frac{n}{\tilde{a}_d \pi_0^2}\tr A^{2d+1}$ only depend on the
spectrum of $\Gamma$, hence so does $\overline{k_d}$ by the
lemma. Therefore, if $\Gamma$ is cospectral with a
distance-regular graph $\Gamma'$, then the average
$\overline{k_d}$ must equal $p_d(k)$, because it does so for
$\Gamma'$. Because the spectrum of a graph determines whether
it is regular and connected, and determines its odd girth, we
hence have:

\begin{cor} {\rm (Huang and Liu~\cite{HL99})} Any graph cospectral with a
generalized odd graph, is a generalized odd graph.
\end{cor}

\section{The odd-girth characterization}\label{odd-girth}

Now let us show that $\overline{k_d} = p_d(k)$ for a connected
regular graph $\Gamma$ having $d+1$ distinct eigenvalues and
finite odd-girth at least $2d+1$. To do this, we need some
basic properties of the predistance polynomials; see also
\cite{Fiolexcess}. First, $\langle p, q\rangle =\frac{1}{n}\tr
(p(A)q(A))$ defines an inner product (determined by the
spectrum of $\Gamma$) on the space of polynomials modulo the
minimal polynomial of $\Gamma$. Using this inner product, one
can find an orthogonal system of so-called predistance
polynomials $p_i, i=0,1,\dots,d$, where $p_i$ has degree $i$
and is normalized such that $\langle p_i, p_i\rangle=p_i(k)
\neq 0$. The predistance polynomials resemble the distance
polynomials of a distance-regular graph; they also satisfy a
three-term recurrence:
$$xp_i=\beta_{i-1}p_{i-1}+\alpha_ip_i+\gamma_{i+1}p_{i+1},\quad
i=0,1,\dots,d,$$ where we let $\beta_{-1}=0$ and $\gamma_{d+1}p_{d+1}=0$
(the latter we may consider as a multiple of the minimal polynomial).
A final property of these polynomials is that $\sum_{i=0}^d p_i$ equals the
Hoffman polynomial $H$.
This implies that the leading coefficient of $p_d$ equals $\frac n{\pi_0}$
(the same as that of $H$).

For the graph $\Gamma$ under consideration, specific properties
hold. It is easy to show by induction that $\alpha_i=0$ for
$i<d$ and that $p_i$ is an even or odd polynomial depending on
whether $i$ is even or odd, for all $i \le d$. Indeed, it is
clear that $p_0=1$ is even and $p_1=x$ is odd, and hence that
$\alpha_0=0$. Now suppose that $\alpha_i=0$ for $i<j<d$ and
that $p_i$ is even or odd (depending on $i$) for $i \le j$.
Then the three-term recurrence implies that
$\alpha_jp_j(k)=\langle xp_j,p_j \rangle=\frac{1}{n}\tr
(Ap_j(A)^2)=0$ because $xp_j^2$ is an odd polynomial of degree
at most $2d-1$. Hence $\alpha_j=0$ and then it follows from the
recurrence that $p_{j+1}$ is even or odd, which finishes the
inductive argument.

What we shall use now is that $xp_d^2$ is an odd polynomial. It
follows that
$$
\alpha_dp_d(k)=\langle xp_d, p_d\rangle=\frac{1}{n}\tr (Ap_d(A)^2) =
\frac{n}{\pi_0^2}\tr A^{2d+1}.
$$
Thus, we have almost shown that this expression for $p_d(k)$
and the one for $\overline{k_d}$ in the lemma are the same;
what remains is to show that $\alpha_d=\tilde{a}_d$. Therefore,
consider again vertices $u$ and $v$ at distance $d$. Then
$$
\alpha_d=\alpha_d(H(A))_{uv}=\alpha_d(p_d(A))_{uv}=(Ap_d(A))_{uv} =
\frac{n}{\pi_0}(A^{d+1})_{uv}=\tilde{a}_d.$$
where the second last step follows because $xp_d$ is odd or
even, and therefore has no term of degree $d$.
Thus, $\overline{k_d}=p_d(k)$ and by the spectral excess theorem we
derive that $\Gamma$ is distance-regular, which finishes the
proof of our result.
\begin{thm}
Let $\Gamma$ be a connected regular graph with $d+1$ distinct
eigenvalues and finite odd-girth at least $2d+1$.
Then $\Gamma$ is a distance-regular generalized odd graph.
\end{thm}
\noindent It is unclear whether we can drop the regularity
condition on $\Gamma$, or in other words, whether there exist
nonregular graphs with $d+1$ distinct eigenvalues and odd-girth
$2d+1$. For nonregular graphs it matters what matrix we
consider (adjacency, Laplacian, etc.). However, for $d=2$ we
know the following:
\begin{prp}
For the adjacency matrix, as well as for the Laplacian matrix, a
connected graph with odd-girth five and three distinct eigenvalues
is regular (and hence distance-regular).
\end{prp}
\noindent {\bf Proof.} For the adjacency matrix $A$ we consider
the minimal polynomial $m$. Suppose $\lambda_0 > \lambda_1 >
\lambda_2$ are the distinct eigenvalues of $A$. The diagonal of
$m(A)=O$ gives that
$(\lambda_0+\lambda_1+\lambda_2)k_u=-\lambda_0\lambda_1\lambda_2$,
where $k_u$ is the valency of vertex $u$. In case
$\lambda_0+\lambda_1+\lambda_2 = \lambda_0\lambda_1\lambda_2 =
0$, it follows that $\lambda_0=-\lambda_2$ and $\lambda_1=0$,
so the graph would be bipartite, which is false. Thus $k_u$ is
constant.

For a graph whose Laplacian matrix has three distinct
eigenvalues it is known that the number $\overline{\mu}$ of
common nonneighbors of two adjacent vertices is constant
(see~\cite{DH98}). Since there are no triangles, it follows
that if $u$ and $v$ are adjacent, then
$k_u+k_v=n-\overline{\mu}$. This implies that any two vertices
at distance two have the same valency. The graph is connected
with at least one odd cycle, hence there exists a walk of even
length between any two vertices $u$ and $v$. Because there are
no triangles, every even vertex on that walk (which includes
$u$ and $v$) has the same valency. \qed

\noindent For the adjacency matrix we also managed to prove
regularity for the analogous cases with four and five distinct
eigenvalues, but we choose not to include the technical
details.\\

\noindent {\bf Acknowledgements} The authors thank the referees
for their useful comments.

\end{document}